\documentclass[12pt]{article}
\usepackage{amssymb}
\usepackage{latexsym,bm}
\usepackage{graphicx}
\usepackage{amsmath}
\usepackage{mathrsfs}

\newcommand{\bea}{\begin{eqnarray*}}
\newcommand{\eea}{\end{eqnarray*}}
\newcommand{\be}{\begin{equation}}
\newcommand{\ee}{\end{equation}}
\newcommand{\ben}{\begin{eqnarray*}}
\newcommand{\een}{\end{eqnarray*}}

\date{}
\begin{document}
\title{Relation between the number of leaves of a tree and its diameter\footnote{E-mail addresses:
{\tt 235711gm@sina.com}(P.Qiao),
{\tt zhan@math.ecnu.edu.cn}(X.Zhan).}}
\author{\hskip -10mm Pu Qiao$^a$, Xingzhi Zhan$^b$\thanks{Corresponding author.}\\
{\hskip -10mm $^a$\small  Department of Mathematics, East China University of Science and Technology, Shanghai 200237, China}\\
{\hskip -10mm $^b$\small Department of Mathematics, East China Normal University, Shanghai 200241, China}}\maketitle
\begin{abstract}
 Let $L(n,d)$ denote the minimum possible number of leaves in a tree of order $n$ and diameter $d.$ In 1975 Lesniak gave the lower bound
 $B(n,d)=\lceil 2(n-1)/d\rceil$ for $L(n,d).$ When $d$ is even, $B(n,d)=L(n,d).$ But when $d$ is odd, $B(n,d)$ is smaller than $L(n,d)$ in general.
 For example, $B(21,3)=14$ while $L(21,3)=19.$ We prove that for $d\ge 2,$
 $$
 L(n,d)=\begin{cases}\left\lceil \frac{2(n-1)}{d}\right\rceil\quad {\rm if}\,\,\,d\,\,\,{\rm is}\,\,\,{\rm even};\\
                    \left\lceil \frac{2(n-2)}{d-1}\right\rceil\quad {\rm if}\,\,\,d\,\,\,{\rm is}\,\,\,{\rm odd}.
                    \end{cases}
 $$
 The converse problem is also considered. Let $D(n,f)$ be the minimum possible diameter of a tree of order $n$ with exactly $f$ leaves.
 We prove that
 $$
 D(n,f)=\begin{cases} 2 \quad\quad\quad\,\,{\rm if}\,\,\,n=f+1;\\
                     2k+1\quad {\rm if}\,\,\,n=kf+2;\\
                     2k+2\quad {\rm if}\,\,\,kf+3\le n\le (k+1)f+1.
                    \end{cases}
 $$
\end{abstract}

\vskip 5mm

A {\it leaf} in a graph is a vertex of degree $1.$ For a real number $r,$ $\lfloor r\rfloor$ denotes the largest integer less than or equal to $r,$
and $\lceil r\rceil$ denotes the least integer larger than or equal to $r.$
Let $L(n,d)$ denote the minimum possible number of leaves in a tree of order $n$ and diameter $d.$
In 1975 Lesniak [1, Theorem 2 on p.285] gave the lower bound  $B(n,d)=\lceil 2(n-1)/d\rceil$ for $L(n,d).$ When $d$ is even, $B(n,d)=L(n,d).$ But when $d$ is odd,
$B(n,d)$ is smaller than $L(n,d)$ in general.  For example, $B(21,3)=14$ while $L(21,3)=19.$

In this note we first determine $L(n,d).$ We use an idea different from that in [1]. The proof also makes it clear why $L(n,d)$
has such an expression. We then determine the minimum possible diameter of a tree with given order and number of leaves.

We  make the necessary preparation. For terminology and notation we follow the books [3] and [2].
We denote by $V(G)$ the vertex set of a graph $G$ and by $d(u,v)$ the distance between two vertices $u$ and $v.$
For vertices $x$ and $y,$ an {\it ($x,y$)-path} is a path with end vertices $x$ and $y.$
We denote by ${\rm deg}(v)$ the degree of a vertex $v.$

Let $P$ be a path in a tree $T$ and we call $P$ the {\it stem} of $T.$ For every vertex $x\in V(T),$ there is a unique ($x,y$)-path $Q$ such that $V(Q)\cap V(P)=\{y\}.$ We say that {\it $x$ originates from $y.$}  Note that by definition, a vertex on the stem originates from itself. A {\it diametral path} of a tree $T$ is a path of length
equal to the diameter of $T.$

A {\it spider} is a tree with at most one vertex of degree larger than $2$ and this vertex is called the {\it branch vertex.} If  no vertex has degree larger than $2,$ then any vertex may be specified as the branch vertex. Thus, a spider is a subdivision of a star. A {\it leg} of a spider is a path from the branch
vertex to a leaf.

We will need the following lemma.

{\bf Lemma 1.} [2, p.63] {\it A path $P=v_0v_1v_2...v_k$ in a tree is a diametral path if and only if for every vertex $x,$
$$
d(x,v_i)\le {\rm min}\{i,\, k-i\}
$$
where $x$ originates from $v_i$ with $P$ as the stem.}

The case $d=1$ for $L(n,d)$ is trivial, since the only tree of diameter $1$ is $K_2$ which has two leaves. Thus it suffices to consider the case $d\ge 2.$

{\bf Theorem 2.} {\it Let $L(n,d)$ denote the minimum possible number of leaves in a tree of order $n$ and diameter $d$ with $d\ge 2.$
Then
$$
 L(n,d)=\begin{cases}\left\lceil \frac{2(n-1)}{d}\right\rceil\quad {\rm if}\,\,\,d\,\,\,{\rm is}\,\,\,{\rm even};\\
                    \left\lceil \frac{2(n-2)}{d-1}\right\rceil\quad {\rm if}\,\,\,d\,\,\,{\rm is}\,\,\,{\rm odd}.
                    \end{cases}
$$}

{\bf Proof.} The idea is to show that for any tree $T,$ there is a corresponding spider with the same order, diameter and number of leaves as $T.$
Hence, to determine $L(n,d)$ it suffices to consider spiders.

If $d=n-1,$ then the tree must be a path which has two leaves. In this case the formula for $L(n,d)$ is true. 
Note also that a path is a spider. 
Next we assume $d\le n-2.$

Let $T$ be a tree of order $n$ and diameter $d.$  Choose a diametral path $P=v_0v_1v_2...v_d$ as the stem.
Suppose that $x$ is a leaf of $T$ outside $P$ originating from $y.$ There is a unique $(x,y)$-path $Q.$ Since $P$ is a diametral path,
$y\neq v_0,\,v_d.$  Hence ${\rm deg}(y)\ge 3.$ We define the {\it first big vertex} of $x,$ denoted by $b(x),$ to be the first vertex of degree at least
$3$ from $x$ to $y$ on $Q.$

Denote $c={\lfloor d/2\rfloor}.$ Then $c=d/2$ if $d$ is even and $c=(d-1)/2$ if $d$ is odd. Let $z=v_c.$
If $T$ has a leaf $u$ outside $P$ with $b(u)\neq z,$ let $w$ be the neighbor of $b(u)$ on the $(b(u),u)$-path. Since $T$ is  a tree, $w$ and $z$
are not adjacent. We delete the edge $wb(u)$ and add the edge $wz$ to obtain a new tree $T_1.$ Since ${\rm min}\{i,\, d-i\}\le {\rm min}\{c,\, d-c\}$
for any $0\le i\le d,$ by Lemma 1 we deduce that $P$ remains a diametral path of $T_1.$ Clearly $T_1$ and $T$ have the same set of leaves.
Hence $T_1$ and $T$ have the same order, diameter and number of leaves. We still designate $P$ as the stem of $T_1$. If $T_1$ has a leaf outside $P$
whose first big vertex is not $z,$ perform the above operation on $T_1$ to obtain a tree $T_2.$  Repeating this operation in the resulting trees
successively finitely many times, we obtain a tree in which every leaf outside $P$ originates from $z$ and with $z$ as its first big vertex.
Such a tree is a spider. An example of the above transformations is depicted in Figure 1.
\vskip 3mm
\par
 \centerline{\includegraphics[width=4in]{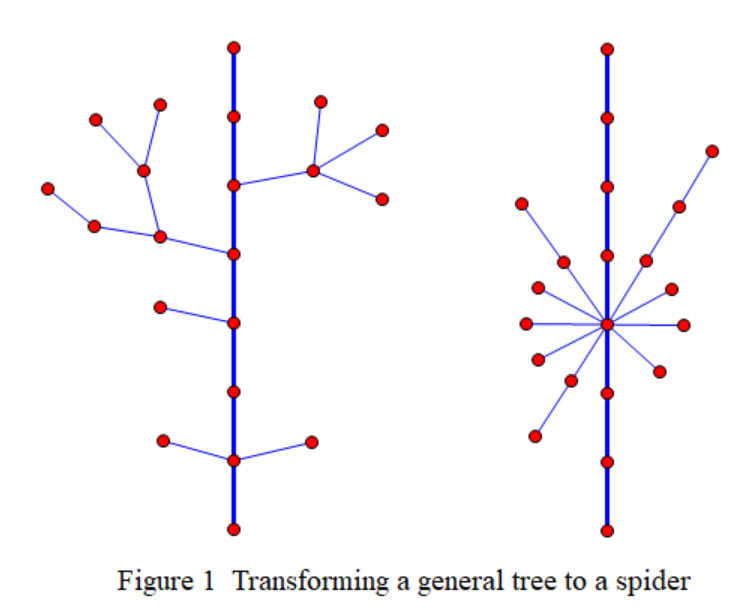}}
\par
The above analysis shows that $L(n,d)$ can be attained at a spider $S$ with a diametral path $P=v_0v_1v_2...v_d$
where $z=v_c$ is the branch vertex. Clearly the number of leaves in $S$ is equal to the number of legs of $S.$ To make the number of legs as small as possible, we need to make each leg as long as possible. Since the diameter of $S$ is $d,$ except the leg $v_cv_{c+1}\ldots v_d$ when $d$ is odd, every other leg has length
at most $c.$  Thus the minimum possible number of legs of such a spider is $\lceil (n-1)/c\rceil$ when $d$ is even and is
$\lceil (n-2)/c\rceil$ when $d$ is odd. This completes the proof. $\Box$

Next we consider the converse problem: Determine the minimum possible diameter of a tree of order $n$ with exactly $f$ leaves.
It suffices to treat the case when $n\ge f+1,$ since $K_2$ is the only tree with $n\le f.$

{\bf Theorem 3.} {\it Let $D(n,f)$ be the minimum possible diameter of a tree of order $n$ with exactly $f$ leaves. Then
$$
 D(n,f)=\begin{cases} 2 \quad\quad\quad\,\, {\rm if}\,\,\,n=f+1;\\
                     2k+1\quad {\rm if}\,\,\,n=kf+2;\\
                     2k+2\quad {\rm if}\,\,\,kf+3\le n\le (k+1)f+1.
                    \end{cases}
 $$}

{\bf Proof.} In the proof of Theorem 2, we showed that for any tree $T,$ there is a corresponding spider with the same order, diameter and number of leaves as $T.$ Thus, it suffices to consider spiders. Note that the number of leaves of a spider is equal to its number of legs, which is also true for the case
when the spider is a path (corresponding to $f=2$) if we take a central vertex of the path as its branch vertex. Let $S$ be a spider of order $n$ with exactly
$f$ legs whose lengths are $x_1\ge x_2\ge\cdots\ge x_f$ arranged in nonincreasing order. Then the diameter of $S$ is $x_1+x_2.$ Hence our problem
is equivalent to minimizing $x_1+x_2$ under the constraint
$$
x_1+x_2+x_3+\cdots+x_f=n-1 \eqno (1)
$$
where $x_1\ge x_2\ge\cdots\ge x_f$ are positive integers.

If $n=f+1,$ then (1) becomes $x_1+x_2+x_3+\cdots+x_f=f,$ which has the only solution $x_1=x_2=x_3=\cdots=x_f=1.$ Hence $x_1+x_2=2.$

Let $n=kf+2.$ If $x_1+x_2\le 2k,$ then $x_2\le k$ and consequently $x_i\le k$ for each $i=3,\ldots,f.$ It follows that
$$
x_1+x_2+x_3+\cdots+x_f\le (x_1+x_2)+(f-2)k\le 2k+(f-2)k=fk=n-2,
$$
contradicting (1). This shows that $D(n,f)\ge 2k+1.$ On the other hand, the values $x_1=k+1,$ $x_2=\cdots=x_f=k$ satisfy (1) and $x_1+x_2=2k+1.$
Hence $D(n,f)=2k+1.$

Now consider the third case $kf+3\le n\le (k+1)f+1.$ We have $kf+2\le n-1\le kf+f.$ Thus there exists an integer $r$ with $2\le r\le f$
such that $n-1=kf+r.$ We first show $D(n,f)\ge 2k+2.$
If $x_1+x_2\le 2k+1,$ then $x_2\le k$ and consequently each $x_i\le k$ for $i=3,\ldots,f.$
It follows that
\begin{eqnarray*}
x_1+x_2+x_3+\cdots+x_f&\le& (x_1+x_2)+(f-2)k\\
                      &\le& 2k+1+(f-2)k\\
                      &=& fk+1\\
                      &<& fk+r=n-1,
\end{eqnarray*}
contradicting (1). Hence $D(n,f)\ge 2k+2.$ On the other hand, the values $x_1=x_2=\cdots=x_r=k+1$ and $x_{r+1}=\cdots=x_f=k$ satisfy
(1) and $x_1+x_2=2k+2,$ which shows $D(n,f)=2k+2.$ This completes the proof. $\Box$

Finally we remark that the maximum problem corresponding to Theorem 2 or Theorem 3 is trivial. The maximum possible number of leaves in a tree of order $n$ and diameter $d$ is $n-d+1$ and the maximum possible diameter of a tree of order $n$ with exactly $f$ leaves is $n-f+1.$

\vskip 5mm
{\bf Acknowledgement.} The authors were supported by the NSFC grants 11671148 and 11771148 and Science and Technology Commission of Shanghai Municipality (STCSM) grant 18dz2271000.

\end{document}